\documentclass[12pt,reqno]{amsart}

\usepackage{amscd,amsmath,amssymb}
\usepackage[mathscr]{eucal}
\usepackage[all]{xy}
\usepackage{units}
\usepackage{blkarray}
\usepackage{hyperref}
\usepackage{xcolor}
\hypersetup{colorlinks=true,urlcolor=black,linkcolor=black,citecolor=black}

\makeatletter \@addtoreset{equation}{section}\makeatother

\newtheorem{theorem}{Theorem}[section]
\newtheorem{lemma}[theorem]{Lemma}
\newtheorem{proposition}[theorem]{Proposition}

\newtheorem{definition}[theorem]{Definition}
\newtheorem{remark}[theorem]{Remark}

\newcommand{\Z}{{\mathbb{Z}}}
\newcommand{\C}{{\mathbb{C}}}
\newcommand{\ZZ}{{\mathbb{Z}/2}}

\newcommand{\ev}{{\mathrm{ev}}}
\newcommand{\od}{{\mathrm{od}}}

\newcommand{\Hom}{{\mathrm{Hom}}}
\newcommand{\End}{{\mathrm{End}}}

\newcommand{\aX}{{X^{\mathtt{h}}}}
\newcommand{\MMFF}{{\mathscr{MF}}}
\newcommand{\MF}{{\mathscr{MF}(f)}}
\newcommand{\HMF}{{\mathfrak{MF}(f)}}
\newcommand{\cA}{{\mathscr{A}}}
\newcommand{\cT}{{\mathfrak{A}}}
\newcommand{\hMF}{{\mathscr{MF}^{\tt h}(f)}}
\newcommand{\dMF}{{\mathscr{MF}^{\mathtt{D}}(f)}}
\newcommand{\cMF}{{\mathscr{MF}_{\tt c}^{\mathtt{D}}(f)}}

\newcommand{\fC}{{\mathsf{C}}}
\newcommand{\fHH}{{\mathsf{HH}}}
\newcommand{\cycH}{{\mathsf{H}}}
\newcommand{\bbf}{{{g}}}
\newcommand{\bbn}{{{0}}}

\newcommand{\cH}{\mathcal{A}}
\newcommand{\rH}{{\mathrm{H}}}
\newcommand{\fH}{{\mathrm{H}}}
\newcommand{\cV}{{\mathscr{H}}}
\newcommand{\tV}{{\cV}}

\newcommand{\ffH}{{{\fH}}}
\newcommand{\fU}{{{U}}}
\newcommand{\cE}{{\mathscr{E}}}
\newcommand{\cD}{{\mathscr{D}}}

\newcommand{\id}{\mathbf{1}}
\newcommand{\pr}{{\pi}}
\newcommand{\ho}{\widehat{\cV}}
\newcommand{\bdp}{\bar{\boldsymbol{\delta}}}
\newcommand{\bbd}{{\boldsymbol{\delta}}}
\newcommand{\bp}{\bar{\partial}}
\newcommand{\bbp}{{\boldsymbol{\bar{\partial}}}}
\newcommand{\Fh}{{{\rm I}^{\sf h}}}
\newcommand{\Fd}{{{\rm I}^{\tt D}}}
\newcommand{\Fc}{{{\rm I}_{\tt c}}}
\newcommand{\Pc}{{\boldsymbol{\pi}}}

\newcommand{\ioh}{{i^{\sf h}}}
\newcommand{\iod}{{i^{\tt D}}}
\newcommand{\ioc}{{i_{\tt c}}}
\newcommand{\str}{{{\sf str}}}
\newcommand{\tr}{{\Theta_{\tt c}}}
\newcommand{\trr}{{\Theta}}
\newcommand{\shi}{{{{\mathit{s}}}}}

\newcommand{\DD}{{\mathrm{D}}}

\allowdisplaybreaks

\title{Calabi-Yau structures on categories of matrix factorizations}
\author{Dmytro Shklyarov\\}

\email{dmytro.shklyarov@mathematik.tu-chemnitz.de}
\dedicatory{\normalsize To Yan Soibelman on the occasion of his 60th anniversary\\}

\begin{document}

\vspace*{-1.0cm}

\begin{abstract} 
We write out explicit proper Calabi-Yau structures, i.~e. non-degenerate cyclic cocycles on the differential graded categories of matrix factorizations of regular functions with isolated critical points.  The formulas involve the Kapustin-Li trace and its ``higher corrections''.
\end{abstract}

\maketitle
\vspace*{-1.0cm}

\baselineskip 1.6pc
\setcounter{tocdepth}{1}
\tableofcontents

\section{Introduction} In the 2000's, following a proposal by Kontsevich, physicists established a precise link between D-branes in B-type topological Landau-Ginzburg models and  matrix factorizations of the corresponding Landau-Ginzburg superpotentials \cite{KL1,BHLS,La4}. Even though the matrix factorizations had, by that time, already been an active area of research in mathematics, the physics interpretation very soon yielded new insights into the subject. One of the first examples of that was a simple universal formula, discovered in \cite{KL2} (``the Kapustin-Li formula''), for Calabi-Yau structures on categories of matrix factorizations. The main result of the present work is a refinement of this formula. In the rest of this introduction we recall what the original formula looks like, explain why (and in what sense) it needs to be refined, and outline our approach to the problem.

Let $f$ be a regular function (a ``superpotential'') on an open affine subset $X\subset\C^n$.  Assume that the only critical value of $f$ is 0  and that all the critical points of $f$ are  isolated. A matrix factorization of $f$ is a $\ZZ$-graded trivial bundle on $X$ endowed with an odd endomorphism $\DD$ satisfying $\DD^2=f\cdot{\rm id}$. Such factorizations can be organized into a differential $\ZZ$-graded category, $\MF$, as follows: the $\ZZ$-graded space of morphisms between two matrix factorizations comprises all homomorphisms between the underlying $\ZZ$-graded bundles and the differential on this space is the super-commutator with the corresponding $\DD$-operators. 

The Kapustin-Li formula describes a Calabi-Yau structure on the $\ZZ$-graded homotopy category $\HMF$ of $\MF$. Let us explain what the term  ``Calabi-Yau structure'' really means in this context. Given a $\ZZ$-graded Hom-finite category $\cT$ (all our categories are $\C$-linear and small), a Calabi-Yau structure of degree (or parity) $d\in\ZZ$ on  $\cT$ is a non-degenerate degree $d$ trace on $\cT$, i.~e. linear maps 
$
\theta: \End^d_\cT(X)\to\C,
$
for all $X\in\cT$, such that the induced pairings
$ 
\Hom^*_\cT(X,Y)\otimes\Hom^{d-*}_\cT(Y,X)\to\C
$
are non-degenerate and graded-symmetric for all $X,Y$.

The category $\HMF$ is known to be Hom-finite (cf. Proposition \ref{mfpr}). The explicit Calabi-Yau structure (of the same parity as $n=\dim\,X$) on $\HMF$ found in \cite{KL2} is given by the formula
\begin{equation}\label{KLtr}
\theta_{\rm KL}: \End^n_{\HMF}(\DD)\to\C,\quad \theta_{\rm KL}(\Phi):=\frac1{n!}\sum_{x}\mathrm{Res}_x\left[\frac{\str\left((\partial\DD)^{\wedge n}\,\Phi\right)}{\partial_1f\ldots\,\partial_nf}\right].
\end{equation}
Here  the summation is over the  critical points of $f$, $\mathrm{Res}_x$ stands for the local residue at $x$, $\str$ for the super-trace, $\partial$ for the holomorphic de Rham differential,  and $\partial_i:=\frac{\partial}{\partial z_i}$. 
(It should be noted that a mathematical proof of the non-degeneracy of the associated pairing -- called by physicists the open-string topological  metric --  was first found only few years later \cite{M,DM}.)

From the physics perspective, the formula (\ref{KLtr}) solves a concrete problem: it completes the description of the open topological field theory \cite{La1,Mo} underlying the B-twisted Landau-Ginzburg model with superpotential $f$. It is at the same time the starting point for a new project \cite{Ca} -- that of ``upgrading'' the open topological field theory to an open topological {\it conformal} field theory \cite{Co,HLL}, i.~e. finding a minimal $A_\infty$ model 
\[
(\HMF, \mu_2\!=\!\boldsymbol{\cdot}, \mu_3,\mu_4,\ldots)
\] for $\MF$ ($\boldsymbol{\cdot}$ denotes the composition of morphisms in $\HMF$)  that is cyclic with respect to $\theta_{\rm KL}$. The latter means that for any $l$ and any  factorizations $\DD^{(i)}$ and morphisms $\Phi_i\in\Hom^*_{\HMF}(\DD^{(i)},\DD^{(i+1)})$ ($i\in \Z/(l+1)\Z$) the expression $\theta_{\rm KL}(\Phi_{l+1}\cdot\mu_{l}(\Phi_{l}\otimes\ldots\otimes\Phi_1))$ has to be cyclically graded-symmetric.

The cyclicity property is a severe restriction on the minimal model and a generic minimal model does not satisfy it. Moreover, the very existence of a cyclic model is far from obvious. Let $\cA$ be a differential $\ZZ$-graded category whose $\ZZ$-graded homotopy category $\cT$ is Hom-finite (in this case $\cA$ is said to be proper) and carries a degree $d$ Calabi-Yau structure $\theta$. How can one be sure that $\cA$ has a minimal $A_\infty$ model that is cyclic with respect to $\theta$? An exhaustive answer to this question was given in \cite{KS}. The answer involves the broader notion of a Calabi-Yau structure on a proper differential $\ZZ$-graded category and can be formulated as follows: A cyclic model exists provided $\theta$ can be ``lifted'' to a Calabi-Yau structure on $\cA$, i.~e. if there is a functional $\Theta:\cycH^\lambda_d(\cA)\to\C$ on the $d$-th cyclic homology group of $\cA$ whose pullback along the composition of certain canonical maps  (cf. Section \ref{cypc})
\begin{equation}\label{lift}
\bigoplus_{X\in\cT}\End^d_{\cT}(X)\to {\sf HH}_d(\cA)\to\cycH^\lambda_d(\cA)
\end{equation}
(${\sf HH}_*$ is the Hochschild homology) coincides with $\theta$. This result  was proven in \cite{KS} using tools of formal non-commutative symplectic geometry. The same tools may be used  -- at least, in principle -- to actually construct a cyclic minimal $A_\infty$ model of $\cA$ starting from any  lift $\Theta$ of $\theta$ and any, not necessarily cyclic minimal model of $\cA$ \cite{Ca}. 

Thus, we are naturally led to the question: Can $\theta_{\rm KL}$ be lifted to $\cycH^\lambda_n(\MF)$? The answer is known to be yes \cite{DM,Se} but this is not at all straightforward. The naive idea that the same formula (\ref{KLtr}), extended  to $\End^n_{\MF}(\DD)$, gives a Calabi-Yau structure on $\MF$ is easily seen to be wrong since,
in general, 
\[
\theta_{\rm KL}(\Phi''\Phi')\neq(-1)^{|\Phi'||\Phi''|}\theta_{\rm KL}(\Phi'\Phi'')
\]
(here $|\cdot|$ denotes the parity of a morphism).  Instead, one can argue as follows: According to \cite{Se}, $\theta_{\rm KL}$ can be extended to a functional on  $\fHH_n(\MF)$ but this already suffices to claim that an extension to $\cycH^\lambda_n(\MF)$ exists as well -- this follows from the degeneration of the Hochschild-to-cyclic spectral sequence for $\MF$ \cite{Dy,DM}.

Given the above-mentioned potential practical significance of Calabi-Yau structures, another natural question to ask is: Are there explicit formulas, similar to (\ref{KLtr}), for a lift of $\theta_{\rm KL}$ to $\cycH^\lambda_n(\MF)$? It is this question that we answer in the present work. The simplest of  our formulas looks as follows:
\begin{equation}\label{corrKL}
\theta_{\rm KL}(\Phi''\Phi')-(-1)^{|\Phi'||\Phi''|}\theta_{\rm KL}(\Phi'\Phi'')=\widetilde{\theta}(\Phi''\otimes\delta(\Phi'))-(-1)^{|\Phi''|}\widetilde{\theta}(\delta(\Phi'')\otimes\Phi')
\end{equation}
where $\Phi'\in \Hom^*_\MF(\DD',\DD'')$, $\Phi''\in \Hom^*_\MF(\DD'',\DD')$, $\delta$ denotes the differential on the morphisms in $\MF$,  and
\begin{multline*}
\widetilde{\theta}(\Psi''\otimes\Psi'):=
\frac{(-1)^n}{(n+1)!}\sum_{x}\sum_{j,\,k=1}^n(-1)^{(k-1)(|\Psi'|+1)}\\
\mathrm{Res}_x\!\left[\frac{\str\left(\Psi''\,(\partial\DD'')^{\wedge k}\,\Psi'\partial_{j} \DD'\,(\partial\DD')^{\wedge(n-k)}+(-1)^{k-1}\Psi''\,\partial_{j} \DD''\,(\partial\DD'')^{\wedge(k-1)}\,\Psi'\,(\partial\DD')^{\wedge (n-k+1)}\right)}{\partial_1f\ldots(\partial_jf)^2\ldots \partial_nf}\right]
\end{multline*}
The functional $\widetilde{\theta}$ is cyclically graded-symmetric:
\[
\widetilde{\theta}(\Psi''\otimes\Psi')=(-1)^{(|\Psi'|+1)(|\Psi''|+1)}\widetilde{\theta}(\Psi'\otimes\Psi''),
\]
which means that $\theta_{\rm KL}+\widetilde{\theta}$ is a kind of ``infinitesimal lift'' of $\theta_{\rm KL}$ to the cyclic homology of $\MF$. It is still not a Calabi-Yau structure and needs to be corrected by ``higher order'' terms. Our main result -- Theorem \ref{mainresult} -- provides explicit formulas for all such higher corrections to $\theta_{\rm KL}$, and thereby solves the problem of lifting the latter to a Calabi-Yau structure on $\MF$. 

In fact, our result is more general, namely, we construct a family of Calabi-Yau structures on $\MF$ depending on a holomorphic volume form $\Omega$ on $X$; the Kapustin-Li trace and its lift to $\MF$ correspond to the special case $\Omega=dz_1\wedge\ldots\wedge dz_n$. The fact that any volume form gives rise to a Calabi-Yau structure on $\MF$   is not surprising  and is in agreement with \cite[Thm.8.3.4]{Pr} where the volume forms were shown to determine  {\it smooth} Calabi-Yau structures on matrix factorization categories. The notion of a smooth Calabi-Yau structure is dual, in some sense, to the one we are interested in here (the latter is also referred to as a {\it proper} Calabi-Yau structure) and for ``nice'' categories, like $\MF$, the two types of Calabi-Yau structures are known to be in bijection  \cite[Prop.6.10]{GPS}.

Let us quickly summarize the main ideas of our approach. To begin with, it does not rely on any of the previously mentioned ideas and results. In particular, the non-degeneracy of the pairing associated with $\theta_{\rm KL}$ is a consequence of the construction. It is also independent of the original paper \cite{KL2} in the sense that the trace $\theta_{\rm KL}$ is not part of the input data. What we do here is an ``open-string'' generalization of the approach to the  closed-string topological metric in the same setting of Landau-Ginzburg models developed in \cite[Sect.2.2]{LLS}.  In \cite{LLS} the authors construct an explicit quasi-isomorphism 
\begin{equation}\label{pv}
({\rm PV}^*(\aX), \lbrace f,\cdot\rbrace)\to({\rm PV}^*(\aX), \lbrace f,\cdot\rbrace)\otimes(\cE_{\tt c}^{(0,*)}(\aX), \bar{\partial})
\end{equation}
where $\aX$ is the analytification of $X$, ${\rm PV}^*(\aX)$ stands for the space of holomorphic polyvector fields, $\lbrace\cdot,\cdot\rbrace$ is the Schouten-Nijenhuis bracket, and
$\cE_{\tt c}^{(0,*)}(\aX)$ is the space of smooth compactly supported $(0,*)$-forms. Any holomorphic volume form on $\aX$ determines a trace on the right-hand side of (\ref{pv}) whose pullback along the quasi-isomorphism can be written in terms of the classical residue trace. This matches the description of the  closed-string topological metric found in \cite{Va}. We apply the same technique to the matrix factorization categories, namely, we construct an explicit $A_\infty$ quasi-equivalence
\begin{equation}\label{mf}
\MF\to\MF\otimes(\cE_{\tt c}^{(0,*)}(\aX), \bar{\partial}),
\end{equation}
observe that any volume form  gives rise to a  Calabi-Yau structure on the right-hand side, and then pull back this Calabi-Yau structure along the quasi-equivalence.

We would like to close by pointing out that the category in the right-hand side of (\ref{mf}) -- we denote it by $\cMF$ in the text -- seems to be an interesting object in its own right. Being an exact analog of the Dolbeault realization of the bounded derived category of a Calabi-Yau variety, $\cMF$ may provide a streamlined approach to various aspects of the B-twisted  Landau-Ginzburg models. For instance, it should allow one to transfer the constructions of \cite{Co1,Po} to the Landau-Ginzburg setting. We hope to return to this topic in the future.

\medskip

\noindent{\bf Conventions.}  For a $\ZZ$-graded space $V$ the parity of an element $v\in V$   will be denoted by $|v|$, $\shi V$ will stand for $V$ with the reversed $\ZZ$-grading, and $\shi$ by itself will stand for the canonical odd map $V\to\shi V$. If $V=(V,d_V)$ is a $\ZZ$-graded complex then $\shi V$ will also denote the complex $(\shi V,d_{\shi V})$ where $d_{\shi V}(\shi v):= -\shi d_V(v)$.  The Koszul sign rule will always be assumed when working with tensors and $[\cdot,\cdot]$ will denote the super-commutator. We will abbreviate ``differential $\ZZ$-graded'' and ``Calabi-Yau'' to  ``dg'' and ``CY''. 

\medskip

\noindent{\bf Acknowledgements.} I would like to thank Manfred Herbst and Daniel Murfet for introducing me to the Kapustin-Li formula and Wolfgang Lerche, Emanuel Scheidegger and Johannes Walcher for inspiring discussions on matrix factorizations and ``open'' mirror symmetry. Special thanks are due to Nils Carqueville for bringing my attention  to the problem of constructing an ``off-shell'' version of the Kapustin-Li trace.

\section{Preliminaries and the main result}

\subsection{CY structures on proper dg categories}\label{cypc}

In this section $\cA$ stands for a small $\C$-linear dg category. 

Let us first recall the definition of the {\it Hochschild complex} $(\fC_*(\cA),b)$ of $\cA$. 
For $l\geq1$ set
\begin{eqnarray*}
\fC^{\{l\}}_*(\cA):=\bigoplus
\Hom^*_\cA(X^{(l)},X^{(1)})\otimes \shi\Hom^*_\cA(X^{(l-1)},X^{(l)})\otimes\ldots\otimes
\shi\Hom^*_\cA(X^{(1)},X^{(2)})
\end{eqnarray*}
where the sum is over all length $l$ collections $X^{(1)},\ldots, \,X^{(l)}$ of objects of $\cA$  and the $\ZZ$-grading on the left-hand side is the total grading on the tensor products. Then, by definition,  
$
\fC_*(\cA)=\bigoplus_{l\geq1} \fC^{\{l\}}_*(\cA).
$
(For $a_{l}\in\Hom^*_\cA(X^{(l)},X^{(1)})$ and $a_i\in\Hom^*_{\cA}(X^{(i)},X^{(i+1)})$, $i\leq l-1$,  the corresponding element of $\fC^{\{l\}}_*(\cA)$ will be written as $a_{l}[a_{l-1}|\ldots |a_1]$.) The differential $b$ is the sum of two anti-commuting differentials $b(\delta)$ and $b(\mu)$ where
\[
b(\delta)(a_{l}[a_{l-1}|\ldots |a_1])=\delta a_{l}[a_{l-1}|\ldots |a_1]+\sum\limits_{i=1}^{l-1}(-1)^{\epsilon_{i+1}}a_{l}[a_{l-1}|\ldots |\delta a_i|\ldots|a_1]
\]
($\delta$ stands for the differential on the Hom-complexes of $\cA$ and $\epsilon_i:=\sum_{j\geq i}|\shi a_{j}|$)
and 
\begin{multline}\label{bmu}
b(\mu)(a_{l}[a_{l-1}|\ldots |a_1])=
(-1)^{|a_{l}|}a_{l}a_{l-1}[a_{l-2}|\ldots |a_1]-\\-\sum\limits_{i=1}^{l-2}(-1)^{\epsilon_{i+1}}a_{l}[a_{l-1}|\ldots
|a_{i+1}a_{i}|\ldots|a_1]-
(-1)^{|\shi a_1|(\epsilon_2+1)}a_1a_{l}[a_{l-1}|\ldots |a_{2}].
\end{multline}

The definition of the {\it cyclic complex} $(\fC_*^\lambda(\cA),b)$ of $\cA$ involves the  cyclic permutation
\begin{equation}\label{cycper}
\tau: \fC_*(\cA)\to \fC_*(\cA),\quad \tau(a_{l}[a_{l-1}|\ldots |a_1])=(-1)^{|\shi a_{l}|(\epsilon_1-|\shi a_{l}|)}a_{l-1}[a_{l-2}|\ldots |a_1|a_{l}].
\end{equation}
One can easily check that $b(\delta)(1-\tau)=(1-\tau)b(\delta)$ and $b(\mu)(1-\tau)=(1-\tau)b'(\mu)$
where $b'(\mu)$ is the operator on $\fC_*(\cA)$ given by the second line in (\ref{bmu}). It follows that the image ${\rm Im}(1-\tau)\subset \fC_*(\cA)$ is $b$-invariant. Then $(\fC_*^\lambda(\cA),b)$ is defined as the quotient of $(\fC_*(\cA),b)$ by the subcomplex $({\rm Im}(1-\tau),b)$. The cohomology of this complex -- the cyclic homology of $\cA$ -- is denoted by $\cycH_*^\lambda(\cA)$.

Recall \cite[Sect.8.2]{KS} that a dg category $\cA$ is said to be {\it proper}  if  for any pair $X'$ and $X''$ of objects $\dim \rH^*(\Hom_{\cA}(X',X''),\delta)<\infty$.

\begin{definition}{\rm (\cite[Sect.10.2]{KS}) Let $\cA$ be a proper dg category. A {\it degree $d\in\ZZ$ CY structure} on $\cA$ is a functional
$
\trr: \cycH_d^\lambda(\cA)\to\C
$ 
with the property that for any objects $X'$ and $X''$ the pairing
\begin{eqnarray}\label{ndc}
\rH^*(\Hom_{\cA}(X'',X'),\delta)\otimes \rH^{d-*}(\Hom_{\cA}(X',X''),\delta)\to \C,\quad
a''\otimes a'\mapsto \trr(\iota(a''a'))
\end{eqnarray}
is non-degenerate. Here $\iota: \bigoplus_X\rH^*(\End_{\cA}(X),\delta)\to \cycH_*^\lambda(\cA)$ is the map induced by the composition of the  embedding $\bigoplus_X(\End^*_{\cA}(X),\delta)=(\fC^{\{1\}}_*(\cA),b(\delta))\to  (\fC_*(\cA),b)$ and the projection $(\fC_*(\cA),b)\to (\fC_*^\lambda(\cA),b)$.
}
\end{definition}

\begin{remark}{\rm We will also speak of {\it chain-level} CY structures. By a chain-level CY structure on $\cA$ we will understand an even/odd functional
$
\trr: \fC_*(\cA)\to\C
$
that descends to $\cycH_*^\lambda(\cA)$, meaning 
\begin{equation}\label{clcy}
\trr\cdot(1-\tau)=0\quad \text{and}\quad\trr\cdot b=0,
\end{equation}
and induces a CY structure in the above sense. }
\end{remark}

\subsection{Matrix factorizations}\label{smf}

Let $(X,f)$ be as in the Introduction, i.~e.
\begin{itemize}
\item $X$ is an open affine subset of $\C^n={\rm Spec}\,\C[z_1,\ldots,z_n]$;
\item $f$ is a regular function on $X$ whose only critical value is $0\in\C$ and whose critical points are all isolated. (The set of critical points will be denoted by $C_f$.)
\end{itemize}
  
Let us reiterate the definition of the dg category $\MF=\MMFF(X,f)$. 
Its objects are block matrices of the form 
\[
\DD=\begin{blockarray}{ccc}
\scriptstyle{k \,{\rm columns}}& \scriptstyle{k \,{\rm columns}}  \\
\begin{block}{[c|c]c}
  0& \DD_{12}  & \scriptstyle{k \,{\rm rows}} \\
\BAhhline{--}
  \DD_{21}& 0 & \scriptstyle{k \,{\rm rows}} \\
  \end{block}
\end{blockarray}
\]
where $\DD_{12}$ and $\DD_{21}$ are $k\times k$ matrices with entries in $\C[X]$ satisfying the conditions 
\begin{eqnarray}\label{mfc}
\DD_{12}\DD_{21}=\DD_{21}\DD_{12}=f\cdot {\id}_k \quad(\Leftrightarrow\,\,   \DD^2=f\cdot{\id}_{2k})
\end{eqnarray}  
(${\id}_k$ stands for the identity $k\times k$ matrix). The space of {\it even} resp. {\it odd} morphisms between two objects
\begin{eqnarray}\label{2obj}
\DD'=\begin{blockarray}{ccc}
\scriptstyle{k} & \scriptstyle{k}  \\
\begin{block}{[c|c]c}
  0& \DD'_{12}  & \scriptstyle{k} \\
\BAhhline{--}
  \DD'_{21}& 0 & \scriptstyle{k} \\
  \end{block}
\end{blockarray}\quad \text{and} \quad
\DD''=\begin{blockarray}{ccc}
\scriptstyle{l} & \scriptstyle{l}  \\
\begin{block}{[c|c]c}
  0& \DD''_{12}  &  \scriptstyle{l}\\
\BAhhline{--}
  \DD''_{21}& 0 & \scriptstyle{l} \\
  \end{block}
\end{blockarray}
\end{eqnarray} 
is the space (in fact, $\C[X]$-module) $\Hom^{\ev}_{\MF}(\DD',\DD'')$ resp. $\Hom^{\od}_{\MF}(\DD',\DD'')$ of block matrices of the form
\begin{eqnarray}\label{evodmor}
\Phi=\begin{blockarray}{ccc}
\scriptstyle{k} & \scriptstyle{k}  \\
\begin{block}{[c|c]c}
  \Phi_{11}& 0  & \scriptstyle{l} \\
\BAhhline{--}
  0 & \Phi_{22} & \scriptstyle{l} \\
  \end{block}
\end{blockarray}
\quad \text{resp.}\quad
\Phi=\begin{blockarray}{ccc}
\scriptstyle{k} & \scriptstyle{k}  \\
\begin{block}{[c|c]c}
  0& \Phi_{12}  & \scriptstyle{l} \\
\BAhhline{--}
\Phi_{21} & 0 & \scriptstyle{l} \\
  \end{block}
\end{blockarray}
\end{eqnarray} 
where $\Phi_{ij}$ are arbitrary $k\times l$ matrices with entries in $\C[X]$; the composition of morphisms is the usual matrix multiplication. Finally, the differential $\delta: \Hom^{\ev/\od}_{\MF}(\DD',\DD'')\to \Hom^{\od/\ev}_{\MF}(\DD',\DD'')$ is defined by the formula 
\begin{eqnarray}\label{diff}
\delta(\Phi):=\DD''\Phi-(-1)^{|\Phi|}\Phi \DD'.
\end{eqnarray}
(Note that $\delta$ is a morphism of $\C[X]$-modules.)

\begin{proposition}\label{mfpr} The dg category $\MF$ is proper. 
\end{proposition}
\noindent This is a special case of \cite[Thm.8.1.1]{Pr} but can also be easily shown directly. Namely, for any pair $\DD',\DD''$ of matrix factorizations the complex of $\mathcal{O}_X$-modules, associated with the complex $(\Hom^*_{\MF}(\DD',\DD''),\delta)$ of $\C[X]$-modules,  is (Zariski) locally contractible outside of the set of critical points of $f$. To see it, pick a point $x\in X\setminus C_f$ and assume that   $\partial_i f(x)\neq0$ for some $i$. Then it follows from the equality  $\DD''^{2}=f\cdot{\id}$ that in the affine neighborhood $U:=\{\partial_i f\neq0\}$ of $x$ 
\begin{equation}\label{homot1}
\DD''\cdot {h}_{\DD''}+{h}_{\DD''}\cdot \DD''=\id,\quad {h}_{\DD''}:=\frac{{\partial_i \DD''}}{{\partial_i f}}.
\end{equation}
Viewing ${h}_{\DD''}$ as an  element of $\End^\od_{\MMFF(U,f)}(\DD'')$,   (\ref{homot1}) means $\delta({h}_{\DD''})={id}_{\DD''}$. As a consequence, for every $\Phi\in \Hom^*_{\MMFF(U,f)}(\DD',\DD'')$ one has $\Phi=\delta {h}(\Phi)+{h}\delta(\Phi)$ where 
$
{h}(\Phi):={h}_{\DD''}\cdot \Phi.
$
Thus, the sheaf of $\mathcal{O}_X$-modules associated with  the $\C[X]$-module $\rH^*(\Hom_{\MF}(\DD',\DD''),\delta)$ is a coherent  sheaf supported on a finite subset of $X$ which implies the claim.

\subsection{Main theorem} The main result of the present work is

\begin{theorem}\label{mainresult}  For any nowhere vanishing holomorphic top degree form $\Omega$ on the analytic space associated with $X$ the following functional $\trr=\trr_{\Omega}$ defines a chain-level CY structure on $\MF$ (of the same parity as  $n=\dim\,X$): 
\begin{equation*}
\trr=\sum_{x\in C_f} \trr_{x}:\fC_*(\MF)\to \C
\end{equation*}
where
\begin{multline}\label{upsx}
\trr_{x}(\Phi_{l}[\Phi_{l-1}|\ldots|\Phi_1]):=\frac1{(n+l-1)!}\sum\limits_{\substack{k_1+\ldots+k_{l}=n-1\\k_1,\ldots,k_l\geq0}}(-1)^{k_1\epsilon_1+\ldots+k_l\epsilon_{l}}
\sum_{i=1}^n(-1)^{i}\\
\sum_{\substack {r_1+\ldots+r_n=l\\ r_j\geq0\, j\neq i, r_i\geq1}}r_1!\,\ldots \,r_n!\sum\limits_{\left(i^{(1)},\ldots,\, i^{(l)}\right)\in \Lambda_n^l(r_1,\ldots, r_n)}\\
\sum\limits_{\left(j^{(l)}_1,\ldots, j^{(l)}_{k_1},\ldots, j^{(1)}_1,\ldots, j^{(1)}_{k_{l}}\right)\in {S}_n^i}{\rm sgn}\left(j^{(l)}_1,\ldots, j^{(l)}_{k_1},\ldots, j^{(1)}_1,\ldots, j^{(1)}_{k_{l}}\right)\\\mathrm{Res}_x\left[\frac{\str\left(\Phi_{l}\partial_{i^{(l)}} \DD^{(l)}\,\partial_{j^{(l)}_1} \DD^{(l)}\ldots \partial_{j^{(l)}_{k_{l}}} \DD^{(l)}\cdot \ldots\cdot \Phi_{1}\partial_{i^{(1)}} \DD^{(1)}\,\partial_{j^{(1)}_1} \DD^{(1)}\ldots \partial_{j^{(1)}_{k_{1}}}\DD^{(1)}\right) \wedge
 \Omega}{(\partial_1f)^{r_1+1}\ldots(\partial_if)^{r_i}\ldots (\partial_nf)^{r_n+1}}\right].
\end{multline}
In this formula
\begin{itemize}
\item $\{\DD^{(i)}\}_{i=1,\ldots, l}$ are arbitrary matrix factorizations, $\Phi_{l}\in\Hom^{\ev/\od}_{\MF}(D^{(l)},D^{(1)})$ and $\Phi_i\in\Hom^{\ev/\od}_{\MF}(D^{(i)},D^{(i+1)})$, $i\leq l-1$;
\item $\epsilon_i:=\sum_{j\geq i}|\shi \Phi_{j}|$;
\item $\Lambda_n^l(r_1,\ldots, r_n)$ denotes the subset in $\{1,\ldots,n\}^l$ of those multi-indices $(i^{(1)},\ldots,\, i^{(l)})$ that contain precisely $r_1$ copies of 1, $r_2$ copies of 2 etc. 

\item ${S}_n^i$ ($i=1,\ldots, n$) stands for the set of all permutations of $(1,2,\ldots,n)\setminus \{i\}$; given an element $(j_1,\ldots,j_{n-1})\in {S}_n^i$, ${\rm sgn}(j_1,\ldots,j_{n-1})$ denotes the sign of the corresponding permutation;
\item $\mathrm{Res}_x$ is the local residue at $x$ and  $\str$ is the supertrace.
\end{itemize}
\end{theorem}
 
The rest  of the paper is devoted to the proof of this theorem. We leave it as an exercise for the reader to show that in the case $l=1$ and $\Omega=dz_1\wedge\ldots\wedge dz_n$ the above formula reproduces the Kapustin-Li trace (\ref{KLtr}). Combining this observation with the condition $\trr(b(\Phi''[\Phi']))=0$, one can derive the formula (\ref{corrKL}).
As yet another exercise, the reader may try to prove (part of) the theorem for functions of one variable ``by hand''. Namely, in the special case $n=1$ the formula (\ref{upsx}) becomes quite simple:
\[  
\trr_{x}(\Phi_{l}[\Phi_{l-1}|\ldots|\Phi_1])=-\mathrm{Res}_x\left[\frac{\str\left(\Phi_{l}\partial_{z} \DD^{(l)}\ldots\Phi_{1}\partial_{z} \DD^{(1)}\right] \wedge
 \Omega}{(\partial_zf)^{l}}\right],
\]
and the properties (\ref{clcy}) can be checked directly.

\section{Proof of the main result}
\subsection{``Dolbeault''  models for the category of matrix factorizations}

\begin{definition}{\rm Let
$\hMF$, $\dMF$ and $\cMF$ be the dg categories with the same objects as $\MF$ and with morphism complexes defined as follows:
\begin{enumerate}
\item
$
\Hom^*_{\hMF}({\cdot,\cdot}):=\left(\Hom^*_{\MF}({\cdot,\cdot})\otimes_{\C[X]} \cH(\aX), \delta\otimes 1\right)
$\\
\noindent where $\aX$ is the complex manifold (analytic space) associated with $X$ and $\cH(\aX)$ is the algebra of holomorphic functions.  
\item 
$
\Hom^*_{\dMF}({\cdot,\cdot}):=\left(\Hom^*_{\MF}({\cdot,\cdot})\otimes_{\C[X]} \cE^{(0,*)}(\aX), \bdp:=\delta\otimes 1+1\otimes \bp\right)
$\\
\noindent where $(\cE^{(0,*)}(\aX),\bar{\partial})$ is the differential $\ZZ$-graded algebra of smooth $(0,*)$-forms on $\aX$ (the $\ZZ$-grading comes from the natural $\Z$-grading). 
\item
$
\Hom^*_{\cMF}({\cdot,\cdot}):=\left(\Hom^*_{\MF}({\cdot,\cdot})\otimes_{\C[X]} \cE_{\tt c}^{(0,*)}(\aX), \bdp:=\delta\otimes 1+1\otimes \bar{\partial}\right)
$\\
\noindent where $\cE_{\tt c}^{(0,*)}(\aX)\subset\cE^{(0,*)}(\aX)$ denotes the subalgebra of compactly supported forms. 
\end{enumerate}
}
\end{definition}

\begin{remark}{\rm
Strictly speaking, $\cMF$ is not a dg category in the conventional sense since it does not have  identity morphisms. However, Proposition \ref{main_s1} below implies that it is weakly unital.
}
\end{remark}
Let us comment on the structure of the Hom-complexes in these categories.

The case of $\hMF$ is clear: we simply allow arbitrary holomorphic functions as entries of the matrices representing the morphisms. The rest of the structure - namely, the composition and the differential - are the same as in the algebraic case. The structure of the Hom-complexes in $\dMF$ and $\cMF$ is slightly more complicated. It is still convenient to think of the morphisms in $\dMF$ as matrices with  entries in the algebra of $(0,*)$-forms on $\aX$. For example, a generic even morphism in this category is the sum of matrices of the following two types:
\begin{equation}\label{evod}
\begin{bmatrix}
  \Phi^\ev_{11}& 0    \\
0 & \Phi^\ev_{22} 
\end{bmatrix}\quad\text{and}\quad
\begin{bmatrix}
  0 & \Phi^\od_{12}  \\
\Phi^\od_{21} & 0
\end{bmatrix}
\end{equation}
where $\Phi^\ev_{ij}$ resp. $\Phi^\od_{ij}$ are matrices with entries in $\cE^{(0,\ev)}(\aX)$ resp. $\cE^{(0,\od)}(\aX)$. The odd morphisms have a similar structure. One should be careful with this description though. Note  that in these terms the composition of morphism does not always coincide with the matrix multiplication. Similarly, the action of the differential $\delta\otimes 1$ is not always given by the formula (\ref{diff}) and the action of the differential $1\otimes \bp$ is not necessarily given by the componentwise action of $\bp$. This is a consequence of the Koszul rule of signs. To avoid confusion in the future, we will denote the composition of morphisms and  the differentials $\delta\otimes1$ and $1\otimes\bp$ in both $\dMF$ and $\cMF$ by $\circ$,  $\bbd$ and $\bbp$, respectively. (Thus, $\bdp=\bbd+\bbp$.)
For example, the composition of 
\[
\Psi=\begin{bmatrix}
  \Psi^\od_{11}& 0    \\
0 & \Psi^\od_{22} 
\end{bmatrix}\in \Hom^\od_{\dMF}(\DD'',\DD'''), \quad
\Phi=\begin{bmatrix}
  0 & \Phi^\ev_{12}  \\
\Phi^\ev_{21} & 0
\end{bmatrix}\in \Hom^\od_{\dMF}(\DD',\DD'')
\]
is the negative of the matrix product:
\[
\Psi\circ\Phi=-\begin{bmatrix}
0&   \Psi^\od_{11}\Phi^\ev_{12}   \\
 \Psi^\od_{22}\Phi^\ev_{21} & 0
\end{bmatrix}\in\Hom^\ev_{\dMF}(\DD',\DD''').
\]
As another example, for an odd  $\Phi$ as above
\begin{equation}\label{krs}
\bbp(\Phi)=-\begin{bmatrix}
  0 & \bp\Phi^\ev_{12}  \\
\bp\Phi^\ev_{21} & 0
\end{bmatrix}.
\end{equation}

\medskip

\subsection{Equivalences between the models}

One has obvious dg functors
\begin{eqnarray}\label{dgf}
\MF\stackrel{\Fh}\longrightarrow
\hMF\stackrel{\Fd}\longrightarrow\dMF\stackrel{\Fc}\longleftarrow\cMF
\end{eqnarray}
induced by the natural embeddings of Hom-complexes
\[
\Hom^*_{\MF}({\cdot,\cdot})\stackrel{\ioh}\hookrightarrow
\Hom^*_{\hMF}({\cdot,\cdot})\stackrel{\iod}\hookrightarrow\Hom^*_{\dMF}({\cdot,\cdot})\stackrel{\ioc}\hookleftarrow\Hom^*_{\cMF}({\cdot,\cdot}).
\]
\begin{proposition}\label{main_s1}
The dg functor $\Fh$, $\Fd$ and $\Fc$ are quasi-equivalences.
\end{proposition}
\noindent{\bf Proof.} We have to show that the embeddings $\ioh$, $\iod$ and $\ioc$ are quasi-isomorphisms. In what follows $\DD'$ and $\DD''$ are two arbitrary matrix factorizations.

\medskip

\noindent\underline{\it Proof for $\ioh$:} Since $\cH(\aX)$ is flat as a $\C[X]$-module \cite[Sect.A1.2]{N}, it suffices to prove that the natural map 
\[
\rH^{*}(\Hom_{\MF}(\DD',\DD''),\delta)\simeq\rH^{*}(\Hom_{\MF}(\DD',\DD''),\delta)\otimes_{\C[X]}\cH(\aX)
\]
is an isomorphism. As we know from the proof of Proposition \ref{mfpr}, the sheaf of $\mathcal{O}_X$-modules, underlying the $\C[X]$-module $\rH^*(\Hom_{\MF}(\DD',\DD''),\delta)$, is a  coherent  sheaf supported at the points of $C_f$. In particular, the support  is proper and, as a consequence of GAGA, the space of global sections of the sheaf does not change upon the analytification.

\medskip

\noindent\underline{\it Proof for $\iod$:} That 
$
\iod
$
is a quasi-isomorphism is an elementary consequence of the Dolbeault theorem (and the fact that $\aX$ is Stein).

\medskip

\noindent\underline{\it Proof for $\ioc$:} To show that $\ioc$ 
is a quasi-isomorphism we will construct an explicit inverse-up-to-homotopy 
 by mimicking an idea from \cite[Sect.2.2]{LLS}. 

Let us fix a Hermitian metric $\langle\cdot,\cdot\rangle$ on the bundle of $(1,0)$-forms on $\aX$. Then one has the following analog of (\ref{homot1}):
on $\aX\setminus C_f$
\begin{equation}\label{homot4}
\DD''\fH_{\DD''}+\fH_{\DD''}\DD''=\id,\quad \fH_{\DD''}:=\frac{\langle\partial \DD'',\,\partial f\rangle}{||\partial f||^2}.
\end{equation}
Here $||\cdot||$ stands for the norm associated with the metric and $\langle\partial \DD'',\,\partial f\rangle$ is the result of applying $\langle {\cdot},\,\partial f\rangle$ to every entry of $\partial \DD''$ (thus,  $\fH_{\DD''}$ is a matrix of smooth functions on $\aX\setminus C_f$). If we view $\fH_{\DD''}$ as an element of  $\End^\od_{\MMFF^{\sf \DD}(\aX\setminus C_f,f)}(\DD'')$ then (\ref{homot4}) means
$
\bbd(\fH_{\DD''})={id}_{\DD''}.
$
By analogy with the algebraic case, we obtain an odd endomorphism 
\begin{equation}\label{mult1}
\ffH(\Phi):=\fH_{\DD''}\circ\Phi
\end{equation} of $\Hom^*_{\MMFF^{\tt D}(\aX\setminus C_f,f)}(\DD',\DD'')$ with the property
\begin{equation}\label{homot3}
\Phi=[\bbd,\ffH](\Phi).
\end{equation}
Furthermore, note that the operator $[\bbp,\ffH]$ on $\Hom^*_{\MMFF^{\tt D}(\aX\setminus C_f,f)}(\DD',\DD'')$ is nilpotent: 
\begin{eqnarray}\label{mult2}
[\bbp,\ffH](\Phi)=(\bbp\ffH+\ffH\bbp)(\Phi)=\bbp(\fH_{\DD''})\circ\Phi=-\bp\fH_{\DD''}\circ\Phi
\end{eqnarray}
(the last equality is a special case of (\ref{krs})).
Together with (\ref{homot3}) this implies that $[\bdp,\ffH]$ is an invertible operator on $\Hom^*_{\MMFF^{\tt D}(\aX\setminus C_f,f)}(\DD',\DD'')$.  Consider the operator
\begin{equation*}\label{Sop}
\tV:=\ffH \cdot [\bdp,\ffH]^{-1}=\ffH \cdot (1+[\bbp,\ffH])^{-1}.
\end{equation*}
Obviously,  $[\bdp,\tV]$ is the identity operator on $\Hom^*_{\MMFF^{\tt D}(\aX\setminus C_f,f)}(\DD',\DD'')$. As a consequence of (\ref{mult1}) and (\ref{mult2}),
\begin{equation}\label{G}
\tV(\Phi)=\cV_{\DD''}\circ\Phi, \quad \cV_{\DD''}=\sum_{i}\fH_{\DD''}\circ(\bp\fH_{\DD''})^{\circ i}.
\end{equation}
 We would like to ``extend'' $\tV$ the whole of $\aX$. 
Let us fix a smooth function $\varrho$ on $\aX$ such that 
\begin{equation}\label{cutoff}
\varrho|_{\fU_1} = 1,\quad\varrho|_{\aX\setminus \fU_2} = 0
\end{equation} for some relatively
compact open neighborhoods (in the analytic topology) ${\fU_1}\subset\overline{\fU}_1\subset{\fU_2}$ of $C_f$  and set 
\begin{eqnarray*}
\ho:=(1-\varrho)\cdot\tV: \Hom^{\ev/\od}_{\MMFF^{\tt D}(\aX,f)}(\DD',\DD'')\to \Hom^{\od/\ev}_{\MMFF^{\tt D}(\aX,f)}(\DD',\DD'')
\end{eqnarray*}
($\ho$ is well-defined on the whole of $\aX$ because of (\ref{cutoff})). Thanks to (\ref{G})  $\ho$ is the operator of left multiplication with an element of $\End^\od_{\MMFF^{\tt D}(\aX,f)}(\DD'')$:  
\begin{equation}\label{hod}
\ho(\Phi)=\ho_{\DD''}\circ\Phi, \quad \ho_{\DD''}=(1-\varrho)\cdot\cV_{\DD''}=(1-\varrho)\cdot\sum_{i}\fH_{\DD''}\circ(\bp\fH_{\DD''})^{\circ i}.
\end{equation}
In particular, $\ho$ preserves the subspace $\Hom^*_{\MMFF^{\tt D}_{\tt c}(\aX,f)}(\DD',\DD'')$. We have
\begin{eqnarray}\label{calc}
[\bdp,\ho]=[\bdp,(1-\varrho)]\tV +(1-\varrho)[\bdp,\tV]=-\bp\varrho\,\tV+(1-\varrho)= 1-(\varrho+\bp\varrho\,\tV).
\end{eqnarray}
Thus, the morphism of complexes 
\begin{eqnarray*}
\pr:=\varrho+\bp\varrho\cdot\tV: \Hom^*_{\MMFF^{\tt D}(\aX,f)}(\DD',\DD'')\to \Hom^*_{\MMFF^{\tt D}_{\tt c}(\aX,f)}(\DD',\DD'')
\end{eqnarray*}
satisfies $id-\ioc\pr=[\bdp,\ho]$ and  $id-\pr\ioc=[\bdp,\ho]$ which
finishes the proof.
\hfill$\blacksquare$

\subsection{An $A_\infty$ functor}
Unlike $\ioc$, its homotopy inverse $\pr$ that we have just constructed is not compatible with the composition of morphisms, i.~e. does not define a dg functor. Our goal now is to promote $\pr$ to an {\it $A_\infty$ functor} $\Pc:\dMF\to\cMF$. That is, we want to find a collection $\{\pi_l\}_{l\geq2}$ of odd maps
\begin{gather*}
\shi\Hom^*_{\dMF}(\DD^{(l)},\DD^{(l+1)})\otimes\ldots\otimes\shi\Hom^*_{\dMF}(\DD^{(2)},\DD^{(3)}) 
\otimes \shi\Hom^*_{\dMF}(\DD^{(1)},\DD^{(2)})\\
\Big\downarrow{\pi_l}
\\
\Hom^*_{\cMF}(\DD^{(1)},\DD^{(l+1)})
\end{gather*}
(for arbitrary $\DD^{(1)},\ldots, \DD^{(l+1)}$) that together with $\pi_1:=\pr\cdot \shi$ satisfy the relations 
\begin{multline}\label{ainfrel}
\sum_{i=1}^{l-1} (-1)^{\epsilon_{i+1}}\left(\pi_{l-i}(\Phi_l|\ldots|\Phi_{i+1})\circ\pi_i(\Phi_i|\ldots|\Phi_{1})-\pi_{l-1}(\Phi_l|\ldots|\Phi_{i+1}\circ\Phi_i|\ldots|\Phi_{1})\right)
=\\
=\bdp\pi_l(\Phi_l|\ldots|\Phi_{1})-\sum_{i=1}^{l} (-1)^{\epsilon_{i+1}}\pi_{l}(\Phi_l| \ldots|\bdp\Phi_i |\ldots|\Phi_{1})
\end{multline}
where $\Phi_i\in \Hom^*_{\dMF}(\DD^{(i)},\DD^{(i+1)})$,  
$
\Phi_l|\ldots|\Phi_{1}:=\shi\Phi_l\otimes\ldots\otimes\shi\Phi_1
$
($\epsilon_{i+1}$ here and in the rest of the paper has the same meaning as in Section \ref{cypc} and in Theorem \ref{mainresult}).

By (\ref{G})  $\pr$ is the operator of left multiplication with an even element, namely
\begin{equation}\label{pid}
\pr(\Phi)=\pi_{\DD''}\circ\Phi, \quad \pi_{\DD''}:=\varrho+\bp\varrho\cdot\cV_{\DD''}=\varrho+\bp\varrho\cdot\sum_{i}\fH_{\DD''}\circ(\bp\fH_{\DD''})^{\circ i}.
\end{equation}
Recall also the odd elements $\ho_{\DD''}$ defined in (\ref{hod}).

\begin{proposition}\label{ainffun} The maps 
\begin{eqnarray*}\label{ainfty}
\pi_l(\Phi_l|\ldots|\Phi_{1})=\pi_{\DD^{(l+1)}}\circ\Phi_l\circ \ho_{\DD^{(l)}}\circ\Phi_{l-1}\circ\ldots \circ\ho_{\DD^{(2)}}\circ \Phi_1
\end{eqnarray*}
define an $A_\infty$ functor $\Pc:\dMF\to\cMF$.
\end{proposition}
\noindent{\bf Proof.} By the Leibniz rule 
\begin{multline*}
\bdp\pi_l(\Phi_l|\ldots|\Phi_{1})=\bdp\left(\pi_{\DD^{(l+1)}}\circ\Phi_l\circ \ho_{\DD^{(l)}}\circ\ldots \circ\ho_{\DD^{(2)}}\circ \Phi_1\right)=\\
=\sum_i (-1)^{\epsilon_{i+1}} \pi_{\DD^{(l+1)}}\circ \ldots\circ\ho_{\DD^{(i+1)}}\circ \bdp\Phi_{i}\circ\ho_{\DD^{(i)}}\circ \ldots \circ \Phi_1-\\
-\sum_i (-1)^{\epsilon_{i+1}} \pi_{\DD^{(l+1)}}\circ \ldots\circ \Phi_{i+1}\circ\bdp\ho_{\DD^{(i+1)}}\circ\Phi_{i}\circ \ldots \circ \Phi_1
\end{multline*}
and therefore
\begin{multline*}
\bdp\pi_l(\Phi_l|\ldots|\Phi_{1})-\sum_{i} (-1)^{\epsilon_{i+1}}\pi_{l}(\Phi_l| \ldots|\bdp\Phi_i |\ldots|\Phi_{1})=\\
=-\sum_i (-1)^{\epsilon_{i+1}}\pi_{\DD^{(l+1)}}\circ \ldots\circ \Phi_{i+1}\circ\bdp\ho_{\DD^{(i+1)}}\circ\Phi_{i}\circ \ldots \circ \Phi_1.
\end{multline*}
By (\ref{calc})  $\bdp\ho_{\DD^{(i+1)}}=id_{\DD^{(i+1)}}-\pi_{\DD^{(i+1)}}$. Hence
\begin{eqnarray*}
&&-\sum_i (-1)^{\epsilon_{i+1}} \pi_{\DD^{(l+1)}}\circ \ldots\circ \Phi_{i+1}\circ\bdp\ho_{\DD^{(i+1)}}\circ\Phi_{i}\circ \ldots \circ \Phi_1=\\
&&=\sum_i (-1)^{\epsilon_{i+1}} \pi_{\DD^{(l+1)}}\circ \ldots\circ \Phi_{i+1}\circ(\pi_{\DD^{(i+1)}}-id_{\DD^{(i+1)}})\circ\Phi_{i}\circ \ldots \circ \Phi_1=\\
&&=\sum_i (-1)^{\epsilon_{i+1}} (\pi_{\DD^{(l+1)}}\circ \ldots\circ \ho_{\DD^{(i+2)}}\circ \Phi_{i+1})\circ(\pi_{\DD^{(i+1)}}\circ\Phi_{i}\circ \ho_{\DD^{(i)}}\circ \ldots \circ \Phi_1)-\\
&&-\sum_i (-1)^{\epsilon_{i+1}} \pi_{\DD^{(l+1)}}\circ \ldots\circ \ho_{\DD^{(i+2)}}\circ(\Phi_{i+1}\circ\Phi_{i})\circ \ho_{\DD^{(i)}}\circ \ldots \circ \Phi_1=\\
&&=\sum_{i} (-1)^{\epsilon_{i+1}}\left(\pi_{l-i}(\Phi_l|\ldots|\Phi_{i+1})\circ\pi_i(\Phi_i|\ldots|\Phi_{1})-\pi_{l-1}(\Phi_l|\ldots|\Phi_{i+1}\circ\Phi_i|\ldots|\Phi_{1})\right). \qquad\blacksquare
\end{eqnarray*}

\begin{remark}\label{ainf}{\rm By Proposition \ref{main_s1} $\Pc$ is not just an $A_\infty$ functor but an $A_\infty$ quasi-equivalence. By precomposing it with the dg quasi-equivalence  $\Fd\cdot\Fh:\MF\to\dMF$ (i.~e. by restricting $\Pc$ to the dg subcategory $\MF$) one obtains an $A_\infty$ quasi-equivalence $\MF\to\cMF$ which we will still denote by $\Pc=\{\pi_l\}_{l\geq1}$.}
\end{remark}

\subsection{CY structures on $\cMF$}

Any holomorphic volume form $\Omega$ on $\aX$ determines a linear functional 
on $\End^*_{\cMF}(\DD)$ (for every $\DD$), namely
\begin{equation}\label{trdef}
\tr(\Phi):=\int_{\aX} \str(\Phi)\wedge \Omega
\end{equation}
where $\str$ is the $\cE_{\tt c}^{(0,*)}(\aX)$-linear extension of the ordinary supertrace 
\[
\str:\End^*_{\MF}(\DD)\to\C[X], \quad \begin{bmatrix}
  \Phi_{11}& \Phi_{12}    \\
\Phi_{21} & \Phi_{22} 
\end{bmatrix}\mapsto \mathsf{tr}(\Phi_{11})-\mathsf{tr}(\Phi_{22}).
\]
Let us extend the resulting functional on $\fC_*^{\{1\}}(\cMF)$ to the whole of $\fC_*(\cMF)$ by setting $\tr|_{\fC_*^{\{l\}}(\cMF)}=0$ for $l\geq2$.  
\begin{proposition} The extended functional $\tr$ is a chain-level CY structure on $\cMF$ (of the same parity as  $n=\dim\,X$). 
\end{proposition}
\noindent{\bf Proof.} The first equality in (\ref{clcy}) is vacuous in this case, while the second one follows from the following easy-to-check properties of $\tr$: 
for any $\Phi\in \End^*_{\cMF}(\DD)$
\begin{equation}\label{p1}
\tr(\bbd\Phi)=\tr(\bbp\Phi)=0\quad(\Rightarrow\,\, \tr(\bdp\Phi)=0)
\end{equation}
and for any $\Phi\in \Hom^*_{\cMF}(\DD',\DD'')$ and $\Psi\in \Hom^*_{\cMF}(\DD'',\DD')$
\begin{equation}\label{p2}
\tr(\Phi\circ\Psi)=(-1)^{|\Phi||\Psi|}\tr(\Psi\circ\Phi).
\end{equation}

That the pairing 
\begin{equation}\label{indpair}
\Hom^*_{\cMF}(\DD'',\DD')\otimes \Hom^*_{\cMF}(\DD',\DD'')\to \C,\quad
\Psi\otimes \Phi\mapsto \tr(\Psi\circ\Phi)
\end{equation}
induces a non-degenerate pairing on the $\bdp$-cohomology is a consequence of the classical Serre duality  \cite{Ser}. Let us sketch the proof. 

Thanks to Proposition \ref{main_s1}, it suffices to show that  
\[
\Hom^*_{\dMF}(\DD'',\DD')\otimes \Hom^*_{\cMF}(\DD',\DD'')\to \C,\quad
\Psi\otimes \Phi\mapsto \tr(\Psi\circ\Phi)
\]
induces a non-degenerate pairing on the $\bdp$-cohomology. 
The spaces $\cE^{(0,*)}(\aX)$ are   Fr\'echet spaces with respect to the topology
of uniform convergence of the forms together with all derivatives on the
compact subsets of $\aX$ \cite[Sect.3]{Ser}. Let $\cD^{(n,*)}(\aX)$ denote the (strong) dual topological vector space of compactly supported currents on $\aX$ of type $(n,*)$. Then the complex $(\Hom^*_{\MF}(\DD'',\DD')\otimes_{\C[X]} \cE^{(0,*)}(\aX), \bdp)$ can be viewed as a 2-periodic complex of Fr\'echet spaces (obviously, both differentials  $\bbd$ and $\bbp$ are continuous maps) and the dual topological complex (cf. \cite[Def.1.1]{LTL}) can be identified with  \[(\Hom^*_{\MF}(\DD',\DD'')\otimes_{\C[X]} \cD^{(n,*)}(\aX), \bdp^\vee=\bbd^\vee+\bbp^\vee)\] where $\vee$ indicates  the transposed map.
By \cite[Thm.1.5,1.6]{LTL} the induced pairing 
\begin{gather*}
\rH^*\left(\Hom_{\MF}(\DD'',\DD')\otimes_{\C[X]} \cE^{(0,*)}(\aX), \bdp\right)\otimes \rH^*\left(\Hom_{\MF}(\DD',\DD'')\otimes_{\C[X]} \cD^{(n,*)}(\aX), \bdp^\vee\right)
\\
\Big\downarrow
\\
\C
\end{gather*}
is non-degenerate because $\rH^*\left(\Hom_{\MF}(\DD'',\DD')\otimes_{\C[X]} \cE^{(0,*)}(\aX), \bdp\right)$ is finite-dimensional (Propositions \ref{mfpr}, \ref{main_s1}) and hence separated.    

By (\ref{p1}) the form $\Omega$ gives rise to a natural embedding 
\[
\left(\Hom^*_{\MF}(\DD',\DD'')\otimes_{\C[X]}\cE_{\tt c}^{(0,*)}, \bdp\right)[n]\to\left(\Hom^*_{\MF}(\DD',\DD'')\otimes_{\C[X]} \cD^{(n,*)}(\aX), \bdp^\vee\right)
\]
and it remains to explain why this embedding is a quasi-isomorphism. In fact, (\ref{p1}) implies more, namely, that the above embedding is a morphism of the underlying {\it double} complexes. Since both double complexes have bounded antidiagonals, it suffices to show that the morphisms 
\[
\left(\Hom^*_{\MF}(\DD',\DD'')\otimes_{\C[X]} \cE_{\tt c}^{(0,*)}, \bbp\right)[n]\to\left(\Hom^*_{\MF}(\DD',\DD'')\otimes_{\C[X]} \cD^{(n,*)}(\aX), \bbp^\vee\right)
\]
are quasi-isomorphisms which is a special case of  \cite[Thm.1]{Ser}.
\hfill$\blacksquare$

\subsection{CY structures on $\MF$}
 
We want to pull back the above CY structure to $\MF$ using the $A_\infty$ quasi-equivalence $\Pc=\{\pi_l\}_{l\geq1}: \MF\to\cMF$ (cf. Remark \ref{ainf}). 

Let $\widehat{\Pc}$ stand for the even linear map $\fC_*(\MF)\to \fC_*^{\{1\}}(\cMF)$ given by 
\[
\widehat{\Pc}(\Phi_{l}[\Phi_{l-1}|\ldots|\Phi_1]):=\pi_{l}(\Phi_{l}|\Phi_{l-1}|\ldots|\Phi_1).
\]
Let also $N$ denote the endomorphism of $\fC_*(\MF)$ that acts as $\sum_{i=0}^{l-1}\tau^i$ on the tensors of length $l$ (here $\tau$ is the cyclic permutation  (\ref{cycper})).

\begin{proposition} The functional
\begin{equation}\label{cymf}
\trr=\tr\cdot \widehat{\Pc}\cdot N: \fC_*(\MF)\to \C
\end{equation}
is a chain-level CY structure on $\MF$ (of the same parity as  $n=\dim\,X$).
\end{proposition}
\noindent{\bf Proof.} The non-degeneracy of the induced pairings 
\[
\rH^*(\Hom_{\MF}(\DD'',\DD'), \delta)\otimes \rH^*(\Hom_{\MF}(\DD',\DD''), \delta)\to\C
\]
follows immediately from the non-degeneracy of the pairings induced by (\ref{indpair}) and the fact that $\Pc$ is an $A_\infty$ quasi-equivalence. The first property in (\ref{clcy}) is obvious since $N\cdot(1-\tau)=0$.  The second property in (\ref{clcy}) follows from \cite[Lem.2.4]{FLS} but we supply an independent proof.

It is easy to see that $N\cdot b(\delta)=b(\delta)\cdot N$ and $N\cdot b(\mu)=b''(\mu)\cdot N$
where
\begin{multline*}
b''(\mu)(\Phi_{l}[\Phi_{l-1}|\ldots |\Phi_1])=
(-1)^{|\Phi_{l}|}\Phi_{l}\Phi_{l-1}[\Phi_{l-2}|\ldots |\Phi_1]-\\-\sum\limits_{i=1}^{l-2}(-1)^{\epsilon_{i+1}}\Phi_{l}[\Phi_{l-1}|\ldots
|\Phi_{i+1}\Phi_{i}|\ldots|\Phi_1].
\end{multline*}
Therefore, $\trr\cdot b=\tr\cdot \widehat{\Pc}\cdot(b(\delta)+b''(\mu))\cdot N$.
Let us calculate the functional $\tr\cdot \widehat{\Pc}\cdot(b(\delta)+b''(\mu))$ explicitly:
\begin{eqnarray*}
&&\tr\cdot \widehat{\Pc}\cdot(b(\delta)+b''(\mu))(\Phi_{l}[\Phi_{l-1}|\ldots |\Phi_1])=\\
&&=\tr\left(\sum\limits_{i=1}^{l}(-1)^{\epsilon_{i+1}}\pi_{l}(\Phi_{l}|\ldots |\delta\Phi_i|\ldots|\Phi_1)-\sum\limits_{i=1}^{l-1}(-1)^{\epsilon_{i+1}}\pi_{l-1}(\Phi_{l}|\ldots
|\Phi_{i+1}\Phi_{i}|\ldots|\Phi_1\right)=\\
&&\stackrel{\bf (\ref{ainfrel})}=\tr\left(\bdp\pi_{l}(\Phi_{l}|\ldots|\Phi_{1})-\sum_{i=1}^{l-1} (-1)^{\epsilon_{i+1}}\pi_{l-i}(\Phi_{l}|\ldots|\Phi_{i+1})\circ\pi_i(\Phi_i|\ldots|\Phi_{1})\right)=\\
&&\stackrel{\bf(\ref{p1})}=-\sum_{i=1}^{l-1} (-1)^{\epsilon_{i+1}}\tr\left(\pi_{l-i}(\Phi_{l}|\ldots|\Phi_{i+1})\circ\pi_i(\Phi_i|\ldots|\Phi_{1})\right).
\end{eqnarray*}
As a result,
\begin{equation*}
\trr\cdot b=\tr\cdot \widehat{\Pc}\cdot(b(\delta)+b''(\mu))\cdot N=\sum_{i=1}^{l-1}\tr\cdot\mu\cdot(\pi_{l-i}\otimes\pi_i)\cdot N
\end{equation*} 
where $\mu$ denotes the composition of morphisms in $\cMF$ and
\begin{equation}\label{pipi}
(\pi_{l-i}\otimes\pi_i)(\Phi_{l}[\Phi_{l-1}|\ldots |\Phi_1]):=(-1)^{\epsilon_{i+1}+1}\pi_{l-i}(\Phi_{l}|\ldots|\Phi_{i+1})\otimes\pi_i(\Phi_i|\ldots|\Phi_{1}).
\end{equation}
Since $\tau\cdot N=N$, one has 
\[
\sum_{i=1}^{l-1}\tr\cdot\mu\cdot(\pi_{l-i}\otimes\pi_i)\cdot N=\sum_{i=1}^{l-1}\tr\cdot\mu\cdot(\pi_{i}\otimes\pi_{l-i})\cdot \tau^{-i}\cdot N.
\]
Observe that
\begin{eqnarray*}
&&\tr\cdot\mu\cdot(\pi_{i}\otimes\pi_{l-i})\cdot \tau^{-i}(\Phi_{l}[\Phi_{l-1}|\ldots |\Phi_1])=\\
&&=
(-1)^{|\shi \Phi_1|(\epsilon_1-|\shi \Phi_1|)+\ldots+|\shi \Phi_i|(\epsilon_1-|\shi \Phi_i|)}\cdot\tr\cdot\mu\cdot(\pi_{i}\otimes\pi_{l-i})(\Phi_{i}[\Phi_{i-1|}|\ldots |\Phi_1|\Phi_{l}|\ldots|\Phi_{i+1}])=\\
&&=
(-1)^{\epsilon_{i+1}(\epsilon_1-\epsilon_{i+1})}\cdot\tr\cdot\mu\cdot(\pi_{i}\otimes\pi_{l-i})(\Phi_{i}[\Phi_{i-1|}|\ldots |\Phi_1|\Phi_{l}|\ldots|\Phi_{i+1}])=\\
&&\stackrel{\bf (\ref{pipi})}=(-1)^{\epsilon_{i+1}(\epsilon_1-\epsilon_{i+1})}\cdot(-1)^{\epsilon_1-\epsilon_{i+1}+1}\cdot\tr\left(\pi_{i}(\Phi_{i}|\ldots|\Phi_{1})\circ\pi_{l-i}(\Phi_{l}|\ldots|\Phi_{i+1})\right)=\\
&&=(-1)^{\epsilon_{i+1}\epsilon_1+\epsilon_{1}+1}\cdot\tr\left(\pi_{i}(\Phi_{i}|\ldots|\Phi_{1})\circ\pi_{l-i}(\Phi_{l}|\ldots|\Phi_{i+1})\right)=\\
&&\stackrel{\bf (\ref{p2})}=(-1)^{\epsilon_{i+1}\epsilon_1+\epsilon_{1}+1}\cdot(-1)^{(\epsilon_{i+1}+1)(\epsilon_1-\epsilon_{i+1}+1)}\cdot\tr\left(\pi_{l-i}(\Phi_{l}|\ldots|\Phi_{i+1})\circ\pi_{i}(\Phi_{i}|\ldots|\Phi_{1})\right)=\\
&&=(-1)^{\epsilon_{i+1}}\cdot\tr\left(\pi_{l-i}(\Phi_{l}|\ldots|\Phi_{i+1})\circ\pi_{i}(\Phi_{i}|\ldots|\Phi_{1})\right)=\\
&&\stackrel{\bf (\ref{pipi})}=-\tr\cdot\mu\cdot(\pi_{l-i}\otimes\pi_{i})(\Phi_{l}[\Phi_{l-1}|\ldots |\Phi_1]).
\end{eqnarray*}
and therefore
\[
\sum_{i=1}^{l-1}\tr\cdot\mu\cdot(\pi_{l-i}\otimes\pi_i)\cdot N=\sum_{i=1}^{l-1}\tr\cdot\mu\cdot(\pi_{i}\otimes\pi_{l-i})\cdot \tau^{-i}\cdot N=-\sum_{i=1}^{l-1}\tr\cdot\mu\cdot(\pi_{l-i}\otimes\pi_i)\cdot N.
\quad \blacksquare
\]
 
\subsection{Explicit formulas in terms of residues} Apart from the  form $\Omega$, the CY structure (\ref{cymf}) depends on the following data:
\begin{itemize}
\item[(a)] the Hermitian metric $\langle\cdot,\cdot\rangle$ on the bundle of $(1,0)$-forms  that we used to construct the ``homotopies'' in (\ref{homot4});
\item[(b)] the neighborhoods ${\fU_1},{\fU_2}$ of $C_f$  and  
 the smooth function $\varrho$  satisfying the conditions (\ref{cutoff}) for these neighborhoods;
\end{itemize}
The proof of Theorem \ref{mainresult} amounts now to evaluating $\trr$ for a special choice of (a) and (b). The goal of this final section is to perform this calculation. (A similar calculation was carried out in \cite{LLS} -- cf. the proof of Proposition 2.5 therein -- but our case is somewhat more involved.) 

Let us start by expanding the formula (\ref{cymf}) in the general case.  As before, we fix some matrix factorizations $\{\DD^{(i)}\}_{i=1,\ldots, l}$ and morphisms $\Phi_{l}\in\Hom^*_{\MF}(\DD^{(l)},\DD^{(1)})$ and $\Phi_i\in\Hom^*_{\MF}(\DD^{(i)},\DD^{(i+1)})$, $i\leq l-1$. 
We have
\begin{multline*}
\trr(\Phi_{l}[\Phi_{l-1}|\ldots|\Phi_1])=\tr\cdot \widehat{\Pc}\cdot N(\Phi_{l}[\Phi_{l-1}|\ldots|\Phi_1])=\\
=\sum_{i}(-1)^{\epsilon_{i}(\epsilon_1-\epsilon_{i})}\tr\cdot\pi_{l}(\Phi_{i-1}|\ldots |\Phi_1|\Phi_{l}|\ldots|\Phi_{i})=\\
\stackrel{\bf Prop.\ref{ainffun}}=\sum_{i}(-1)^{\epsilon_{i}(\epsilon_1-\epsilon_{i})}\tr(\pi_{\DD^{(i)}}\circ\Phi_{i-1}\circ\ho_{\DD^{(i-1)}}\circ \ldots \circ   \Phi_1\circ \ho_{\DD^{(1)}}\circ \Phi_{l}\circ \ho_{\DD^{(l)}}\circ \ldots\circ \ho_{\DD^{(i+1)}}\circ \Phi_{i})=\\
\stackrel{\bf (\ref{p2})}=\sum_{i=1}^{l}(-1)^{\epsilon_1-\epsilon_{i}}\tr((\Phi_{l}\circ \ho_{\DD^{(l)}})\circ \ldots\circ  (\Phi_{i}\circ \pi_{\DD^{(i)}})\circ\ldots\circ  (\Phi_1\circ \ho_{\DD^{(1)}}))=\\
\stackrel{\bf (\ref{hod}),(\ref{pid})}=\sum_{i=1}^{l}(-1)^{\epsilon_1-\epsilon_{i}}\tr\left((\Phi_{l}\circ\cV_{\DD^{(l)}})\circ \ldots\circ\Phi_i\circ\ldots\circ  (\Phi_1\circ \cV_{\DD^{(1)}}) \varrho(1-\varrho)^{l-1}\right)-\\
-\sum_{i=1}^{l}\tr\left((\Phi_{l}\circ\cV_{\DD^{(l)}})\circ \ldots\circ  (\Phi_1\circ \cV_{\DD^{(1)}}) \bp\varrho(1-\varrho)^{l-1}\right)=\\
\stackrel{\bf (\ref{trdef})}=\sum_{i=1}^{l}(-1)^{\epsilon_1-\epsilon_{i}}\int_{\aX} \str\left((\Phi_{l}\circ\cV_{\DD^{(l)}})\circ \ldots\circ\Phi_i\circ\ldots\circ  (\Phi_1\circ \cV_{\DD^{(1)}})\right) \varrho(1-\varrho)^{l-1}\wedge \Omega-\\
-l\int_{\aX} \str\left((\Phi_{l}\circ\cV_{\DD^{(l)}})\circ \ldots\circ  (\Phi_1\circ \cV_{\DD^{(1)}})\right) \bp\varrho(1-\varrho)^{l-1}\wedge \Omega=\\
=\sum_{i=1}^{l}(-1)^{\epsilon_1-\epsilon_{i}}\int_{\aX} \str\left((\Phi_{l}\circ\cV_{\DD^{(l)}})\circ \ldots\circ\Phi_i\circ\ldots\circ  (\Phi_1\circ \cV_{\DD^{(1)}})\right) \varrho(1-\varrho)^{l-1}\wedge \Omega+\\
+\int_{\aX} \str\left((\Phi_{l}\circ\cV_{\DD^{(l)}})\circ \ldots\circ  (\Phi_1\circ \cV_{\DD^{(1)}})\right) \bp\left((1-\varrho)^{l}\right)\wedge \Omega.
\end{multline*}

From now on, 
\begin{itemize}
\item[(a)] $\langle\sum_ig'_idz_i,\sum_jg''_jdz_j\rangle:=\sum_ig'_i\overline{g''_i}$;
\item[(b)] ${\fU_2}=\cup_{x\in C_f} B_{2r}(x)$ and ${\fU_1}=\cup_{x\in C_f} B_{r}(x)$ where $B_{r}(x)$ stands for the open ball of radius $r$ centered at $x$ (for the standard metric on $\C^n$). We assume that $r$ is small enough, so that  $\overline{B_{2r}(x)}\subset \aX$ and $\overline{B_{2r}(x)}\cap \overline{B_{2r}(y)}=\varnothing$ for two different $x,y\in C_f$.
\end{itemize}
Then
\begin{multline*}
\sum_{i=1}^{l}(-1)^{\epsilon_1-\epsilon_{i}}\int_{\aX} \str\left((\Phi_{l}\circ\cV_{\DD^{(l)}})\circ \ldots\circ\Phi_i\circ\ldots\circ  (\Phi_1\circ \cV_{\DD^{(1)}})\right) \varrho(1-\varrho)^{l-1}\wedge \Omega+\\
+\int_{\aX} \str\left((\Phi_{l}\circ\cV_{\DD^{(l)}})\circ \ldots\circ  (\Phi_1\circ \cV_{\DD^{(1)}})\right) \bp\left((1-\varrho)^{l}\right)\wedge \Omega=\\
=\sum_{i=1}^{l}(-1)^{\epsilon_1-\epsilon_{i}}\int_{\overline{U}_2\setminus U_1} \str\left((\Phi_{l}\circ\cV_{\DD^{(l)}})\circ \ldots\circ\Phi_i\circ\ldots\circ  (\Phi_1\circ \cV_{\DD^{(1)}})\right) \varrho(1-\varrho)^{l-1}\wedge \Omega+\\
+\int_{\overline{U}_2\setminus U_1} \str\left((\Phi_{l}\circ\cV_{\DD^{(l)}})\circ \ldots\circ  (\Phi_1\circ \cV_{\DD^{(1)}})\right) \bp\left((1-\varrho)^{l}\right)\wedge \Omega
\end{multline*}
which by the Stokes theorem equals 
\begin{multline}\label{finfor}
\sum_{i=1}^{l}(-1)^{\epsilon_1-\epsilon_{i}}\int_{\overline{U}_2\setminus U_1} \str\left((\Phi_{l}\circ\cV_{\DD^{(l)}})\circ \ldots\circ\Phi_i\circ\ldots\circ  (\Phi_1\circ \cV_{\DD^{(1)}})\right) \varrho(1-\varrho)^{l-1}\wedge \Omega+\\
+(-1)^n\int_{\overline{U}_2\setminus U_1} \bp\left(\str\left((\Phi_{l}\circ\cV_{\DD^{(l)}})\circ \ldots\circ  (\Phi_1\circ \cV_{\DD^{(1)}})\right) \right)(1-\varrho)^{l}\wedge \Omega+\\
+(-1)^{n-1}\sum_{x\in C_f}\int_{\partial\overline{B_{2r}(x)}-\partial\overline{B_{r}(x)}} \str\left((\Phi_{l}\circ\cV_{\DD^{(l)}})\circ \ldots\circ  (\Phi_1\circ \cV_{\DD^{(1)}})\right) (1-\varrho)^{l}\wedge \Omega.
\end{multline}
 Taking into account (\ref{cutoff}), the last sum in (\ref{finfor}) is equal to  
\begin{eqnarray*}
(-1)^{n-1}\sum_{x\in C_f}\int_{\partial\overline{B_{2r}(x)}} \str\left((\Phi_{l}\circ\cV_{\DD^{(l)}})\circ \ldots\circ  (\Phi_1\circ \cV_{\DD^{(1)}})\right)\wedge \Omega.
\end{eqnarray*}

\begin{lemma} The first and the second lines in (\ref{finfor}) vanish.
\end{lemma}
\noindent{\bf Proof.} It follows from the definition (\ref{G}) of $\cV_{\DD}$ and from the holomorphicity of the $\Phi$'s that the top degree parts of the integrands in the first and the second lines in (\ref{finfor}) are sums of expressions of the form
\begin{eqnarray*}
\str\left(\xi_{l} \circ (\bp\fH_{\DD^{(l)}})^{\circ k_{l}}\circ \ldots\circ\xi_1\circ  (\bp\fH_{\DD^{(1)}})^{\circ k_{1}}\right) \wedge\ldots
\end{eqnarray*}
where $\xi_i$ are some morphisms and $k_1+\ldots+k_{l}=n$ ($k_i\geq0$). Therefore, it is enough to show that the wedge-product of arbitrary $n$ matrix elements of the matrices $\bp\fH_{\DD^{(i)}}$ equals 0. In view of (\ref{homot4}), it suffices to prove that for any collection $\{\omega_1,\ldots,\omega_n\}$ of holomorphic 1-forms on $\aX$
\[
\bp\frac{\langle\omega_1,\partial f\rangle}{||\partial f||^2}\wedge\ldots\wedge \bp\frac{\langle\omega_n,\partial f\rangle}{||\partial f||^2}=0.
\]
Note that the left-hand side of the latter equality is $\cH(\aX)$-multilinear and skew-symmetric in the $\omega$'s. Thus, we only need to prove it for $\omega_i=dz_i$, $i=1,\ldots n$. We have
\begin{multline*}
\bp\frac{\langle dz_1,\partial f\rangle}{||\partial f||^2}\wedge\ldots\wedge \bp\frac{\langle dz_n,\partial f\rangle}{||\partial f||^2}=\bp\frac{\overline{\partial_1 f}}{||\partial f||^2}\wedge\ldots\wedge \bp\frac{\overline{\partial_n f}}{||\partial f||^2}
=\frac1{\prod_i\partial_i f}\bp\frac{|\partial_1 f|^2}{||\partial f||^2}\wedge\ldots\wedge \bp\frac{|\partial_n f|^2}{||\partial f||^2}=\\
=\frac1{\prod_i\partial_i f}\bp\frac{|\partial_1 f|^2}{||\partial f||^2}\wedge\ldots\wedge \bp\frac{|\partial_{n-1} f|^2}{||\partial f||^2}\wedge \bp\left(1-\frac{|\partial_1 f|^2+\ldots+|\partial_{n-1} f|^2}{||\partial f||^2}\right)=0.
\end{multline*}
\hfill$\blacksquare$

\bigskip

The conclusion so far is that 
\begin{eqnarray}\label{fin1}
\trr(\Phi_{l}[\Phi_{l-1}|\ldots|\Phi_1])=(-1)^{n-1}\sum_{x\in C_f}\int_{\partial\overline{B_{2r}(x)}} \str\left(\Phi_{l}\circ\cV_{\DD^{(l)}}\circ \ldots\circ  \Phi_1\circ \cV_{\DD^{(1)}}\right)\wedge \Omega.
\end{eqnarray}
By (\ref{G}) and (\ref{homot4})
\begin{multline*}
(-1)^{n-1}\int_{\partial\overline{B_{2r}(x)}} \str\left(\Phi_{l}\circ\cV_{\DD^{(l)}}\circ \ldots\circ  \Phi_1\circ \cV_{\DD^{(1)}}\right)\wedge \Omega=\\
=(-1)^{n-1}\sum\limits_{\substack{k_1+\ldots+k_{l}=n-1\\k_1,\ldots,k_l\geq0}}\int_{\partial\overline{B_{2r}(x)}}\\ \str\left(\Phi_{l}\circ\fH_{\DD^{(l)}}\circ(\bp\fH_{\DD^{(l)}})^{\circ k_{l}}\circ \ldots\circ  \Phi_1\circ \fH_{\DD^{(1)}}\circ(\bp\fH_{\DD^{(1)}})^{\circ k_1}\right)\wedge \Omega=\\
=(-1)^{n-1}\sum\limits_{\substack{k_1+\ldots+k_{l}=n-1\\k_1,\ldots,k_l\geq0}}\int_{\partial\overline{B_{2r}(x)}} \\
\str\left(\frac{\langle\Phi_l\partial \DD^{(l)},\partial f\rangle}{||\partial f||^2}\circ\left(\bp\frac{\langle\partial \DD^{(l)},\partial f\rangle}{||\partial f||^2}\right)^{\circ k_{l}}\circ \ldots\circ  \frac{\langle \Phi_1\partial \DD^{(1)},\partial f\rangle}{||\partial f||^2}\circ\left(\bp\frac{\langle\partial \DD^{(1)},\partial f\rangle}{||\partial f||^2}\right)^{\circ k_{1}}\right)\wedge \Omega.
\end{multline*}
Setting $\zeta_j:=\frac{\overline{\partial_j f}}{||\partial f||^2}$ and unfolding the definition of $\circ$, the latter equals
\begin{multline*}
-(-1)^{\frac{n(n+1)}2}\sum\limits_{\substack{k_1+\ldots+k_{l}=n-1\\k_1,\ldots,k_l\geq0}}(-1)^{k_1\epsilon_1+\ldots+k_l\epsilon_{l}}\sum\limits_{\left(j^{(1)}_1,\ldots, j^{(1)}_{k_1},\ldots, j^{(l)}_1,\ldots, j^{(l)}_{k_{l}}\right)}\sum\limits_{\left(i^{(1)}, \ldots, i^{(l)}\right)}\\
\int_{\partial\overline{B_{2r}(x)}} \zeta_{i^{(l)}}\ldots\zeta_{i^{(1)}}\bp\zeta_{j^{(l)}_1}\wedge\ldots\wedge\bp\zeta_{j^{(l)}_{k_{l}}}\wedge\ldots
\wedge \bp\zeta_{j^{(1)}_1}\wedge\ldots\wedge\bp\zeta_{j^{(1)}_{k_{1}}}\wedge\\
\wedge\str\left(\Phi_{l}\partial_{i^{(l)}} \DD^{(l)}\,\partial_{j^{(l)}_1} \DD^{(l)}\ldots \partial_{j^{(l)}_{k_{l}}} \DD^{(l)}\cdot \ldots\cdot \Phi_{1}\partial_{i^{(1)}} \DD^{(1)}\,\partial_{j^{(1)}_1} \DD^{(1)}\ldots \partial_{j^{(1)}_{k_{1}}}\DD^{(1)}\right) \wedge
 \Omega
\end{multline*}
where $i^{(s)}$ and $j^{(s)}_r$ run from 1 to $n$. Since $\bp\zeta_{j^{(l)}_1}\wedge\ldots\wedge\bp\zeta_{j^{(l)}_{k_{l}}}\wedge\ldots
\wedge \bp\zeta_{j^{(1)}_1}\wedge\ldots\wedge\bp\zeta_{j^{(1)}_{k_{1}}}$ is non-trivial only for 
$\left(j^{(l)}_1,\ldots, j^{(l)}_{k_1},\ldots, j^{(1)}_1,\ldots, j^{(1)}_{k_{l}}\right)\in {S}_n^i$ (see Theorem \ref{mainresult} for the definition of ${S}_n^i$),  the above expression is equal to
\begin{multline*}
-(-1)^{\frac{n(n+1)}2}\sum\limits_{\substack{k_1+\ldots+k_{l}=n-1\\k_1,\ldots,k_l\geq0}}(-1)^{k_1\epsilon_1+\ldots+k_l\epsilon_{l}}\sum_{i=1}^n \\
\sum\limits_{\left(j^{(l)}_1,\ldots, j^{(l)}_{k_1},\ldots, j^{(1)}_1,\ldots, j^{(1)}_{k_{l}}\right)\in {S}_n^i} {\rm sgn}\left(j^{(l)}_1,\ldots, j^{(l)}_{k_1},\ldots, j^{(1)}_1,\ldots, j^{(1)}_{k_{l}}\right)
\sum\limits_{\left(i^{(1)}, \ldots, i^{(l)}\right)}\\
\int_{\partial\overline{B_{2r}(x)}} \zeta_{i^{(l)}}\ldots\zeta_{i^{(1)}}\bp\zeta_{1}\wedge\ldots\wedge\widehat{\overline{\partial}\zeta_{i}}\wedge\ldots
\wedge \bp\zeta_{n}\wedge\\
\wedge\str\left(\Phi_{l}\partial_{i^{(l)}} \DD^{(l)}\,\partial_{j^{(l)}_1} \DD^{(l)}\ldots \partial_{j^{(l)}_{k_{l}}} \DD^{(l)}\cdot \ldots\cdot \Phi_{1}\partial_{i^{(1)}} \DD^{(1)}\,\partial_{j^{(1)}_1} \DD^{(1)}\ldots \partial_{j^{(1)}_{k_{1}}}\DD^{(1)}\right) \wedge
 \Omega
\end{multline*}
where $\widehat{\,\,\,}$ means the term is omitted. As it is shown in \cite[Ch.5, Sect.1]{GH}, for any $i=1,\ldots,n$
\begin{eqnarray*}
\bp\zeta_{1}\wedge\ldots\wedge\widehat{\overline{\partial}\zeta_{i}}\wedge\ldots
\wedge \bp\zeta_{n}=(-1)^{i-1} \cdot \partial_if \cdot\boldsymbol{\eta}
\end{eqnarray*}
where $\boldsymbol{\eta}$ stands for the $(0,n-1)$-form
\begin{eqnarray*}
\frac1{||\partial f||^{2n}}\cdot \sum_{s=1}^n(-1)^{s-1} \overline{\partial_s f} \,\overline{\partial (\partial_1 f)}\wedge\ldots\wedge\widehat{\overline{\partial (\partial_s f)}}\wedge\ldots\wedge\overline{\partial(\partial_n f)}.
\end{eqnarray*}
Thus, 
\begin{multline}\label{pref}
\trr(\Phi_{l}[\Phi_{l-1}|\ldots|\Phi_1])=-(-1)^{\frac{n(n+1)}2}\sum\limits_{\substack{k_1+\ldots+k_{l}=n-1\\k_1,\ldots,k_l\geq0}}(-1)^{k_1\epsilon_1+\ldots+k_l\epsilon_{l}}\sum_{i=1}^n(-1)^{i-1}\\
\sum\limits_{\left(j^{(l)}_1,\ldots, j^{(l)}_{k_1},\ldots, j^{(1)}_1,\ldots, j^{(1)}_{k_{l}}\right)\in {S}_n^i} {\rm sgn}\left(j^{(l)}_1,\ldots, j^{(l)}_{k_1},\ldots, j^{(1)}_1,\ldots, j^{(1)}_{k_{l}}\right)
\sum\limits_{\left(i^{(1)}, \ldots, i^{(l)}\right)}\\
\int_{\partial\overline{B_{2r}(x)}} \frac{\overline{\partial_{i^{(l)}} f}}{||\partial f||^2}\ldots\frac{\overline{\partial_{i^{(1)}} f}}{||\partial f||^2} \cdot \boldsymbol{\eta}\\
\wedge\partial_if\cdot\str\left(\Phi_{l}\partial_{i^{(l)}} \DD^{(l)}\,\partial_{j^{(l)}_1} \DD^{(l)}\ldots \partial_{j^{(l)}_{k_{l}}} \DD^{(l)}\cdot \ldots\cdot \Phi_{1}\partial_{i^{(1)}} \DD^{(1)}\,\partial_{j^{(1)}_1} \DD^{(1)}\ldots \partial_{j^{(1)}_{k_{1}}}\DD^{(1)}\right) \wedge
 \Omega.
\end{multline}
Finally, using the decomposition 
\[
\{1,\ldots,n\}^l=\coprod_{\substack {r_1+\ldots+r_n=l\\ r_1,\ldots,r_n\geq0}}\Lambda_n^l(r_1,\ldots, r_n).
\]
(see Theorem \ref{mainresult} for the definition of $\Lambda_n^l(r_1,\ldots, r_n)$) and the formula (\ref{GGH}) with $g_i:=\partial_if$, (\ref{pref}) can be written as
\begin{multline*}
-(-1)^{\frac{n(n+1)}2}(2\pi i)^n\sum\limits_{\substack{k_1+\ldots+k_{l}=n-1\\k_1,\ldots,k_l\geq0}}(-1)^{k_1\epsilon_1+\ldots+k_l\epsilon_{l}}\sum_{i=1}^n(-1)^{i-1}\\
\sum\limits_{\left(j^{(l)}_1,\ldots, j^{(l)}_{k_1},\ldots, j^{(1)}_1,\ldots, j^{(1)}_{k_{l}}\right)\in {S}_n^i}{\rm sgn}\left(j^{(l)}_1,\ldots, j^{(l)}_{k_1},\ldots, j^{(1)}_1,\ldots, j^{(1)}_{k_{l}}\right)\\\sum_{\substack {r_1+\ldots+r_n=l\\ r_1,\ldots,r_n\geq0}}\frac{r_1!\,\ldots \,r_n!}{(n+l-1)!}\sum\limits_{\left(i^{(1)},\ldots,\, i^{(l)}\right)\in \Lambda_n^l(r_1,\ldots, r_n)}\\
\mathrm{Res}_x\left[\frac{\str\left(\Phi_{l}\partial_{i^{(l)}} \DD^{(l)}\,\partial_{j^{(l)}_1} \DD^{(l)}\ldots \partial_{j^{(l)}_{k_{l}}} \DD^{(l)}\cdot \ldots\cdot \Phi_{1}\partial_{i^{(1)}} \DD^{(1)}\,\partial_{j^{(1)}_1} \DD^{(1)}\ldots \partial_{j^{(1)}_{k_{1}}}\DD^{(1)}\right) \wedge
 \Omega}{(\partial_1f)^{r_1+1}\ldots(\partial_if)^{r_i}\ldots (\partial_nf)^{r_n+1}}\right]
\end{multline*}  
which equals, up to a constant,  the right-hand side of (\ref{upsx}).

\appendix

\section{A generalization of a formula of Griffiths and Harris}

Let $\bbf=(g_1,\ldots, g_n): (U,\bbn)\to(\C^n,\bbn)$ be a holomorphic map defined in some open neighborhood $U$ of the origin $\bbn\in\C^n$. Suppose $\bbn$ is an {\it isolated} zero of $\bbf$ and choose an open ball $B$ centered at $\bbn$ so that $\overline{B}\subset U$  and $\overline{B}$ contains no other zeros of $\bbf$. The aim of this appendix is to prove
\begin{proposition} For any holomorphic $n$-form $\omega$ in $\overline{B}$ and non-negative integers $r_i$ 
\begin{eqnarray}\label{GGH}
\mathrm{Res}_\bbn \left[\frac{\omega}{g_1^{r_1+1}\ldots\, g_n^{r_n+1}}\right]=\frac1{(2\pi i)^n}\frac{(n-1+r_1+\ldots+r_n)!}{r_1!\,\ldots \,r_n!}\int_{\partial \overline{B}} \,\frac{\overline{g}_1^{r_1}\ldots\, \overline{g}_n^{r_n}}{||\bbf||^{2(r_1+\ldots+r_n)}}\,\boldsymbol{\eta}(\bbf) \wedge \omega\nonumber\\
\end{eqnarray}
where 
$||\bbf||^2:=|g_1|^2+\ldots+|g_n|^2$  and
\begin{eqnarray*}
\boldsymbol{\eta}(\bbf):=\frac1{||\bbf||^{2n}}\cdot \sum_{s=1}^n(-1)^{s-1} \overline{g}_s \overline{\partial g_1}\wedge\ldots\wedge\widehat{\overline{\partial g_s}}\wedge\ldots\wedge\overline{\partial g_n}
\end{eqnarray*}
(as before, $\widehat{\,\,\,}$ means the term is omitted).
\end{proposition}

The case $r_1=\ldots=r_n=0$ of the formula is due to Griffiths and Harris \cite[Ch.5, Sect.1]{GH}. As we will see, this special case, or rather its global version, implies (\ref{GGH}).

Let us fix an $n$-tuple $\varepsilon=(\varepsilon_1,\ldots, \varepsilon_n)$ of positive numbers and denote by 
\[
D_\varepsilon=\{w=(w_1,\ldots,w_n)\in \C^n\,|\,|w_i|<\varepsilon_i\,\,\forall i\}\subset \C^n
\]
the corresponding polydisc. We will assume that the $\varepsilon_i$ are small enough, so that $\overline{D}_\varepsilon\cap \bbf(\partial \overline{B})=\varnothing$. Then the number of solutions,  counted with multiplicities, of the system $\bbf=w$ in $B$ is the same (finite) number for all $w$ in $D_\varepsilon$ \cite[Sect.5.4]{AGV}.   Let $\omega$ be any holomorphic $n$-form in $\overline{B}$. The main idea of the proof of (\ref{GGH}) is to compare two integral representations of  the {\it global} residue
\begin{eqnarray}\label{slr}
\sum_{x\in \bbf^{-1}(w)\cap B} {\rm Res}_x\left(\frac{\omega}{(g_1-w_1)\ldots\, (g_n-w_n)}\right),\quad w\in D_\varepsilon.
\end{eqnarray}
On the one hand, by the global residue theorem \cite[Ch.5, Sect.1]{GH}
\begin{eqnarray}\label{w1}
(\ref{slr})=\frac{(n-1)!}{(2\pi i)^n}\int_{\partial \overline{B}} \boldsymbol{\eta}(\bbf-w) \wedge \omega
\end{eqnarray}
where
\[
\boldsymbol{\eta}(\bbf-w)=\frac{1}{||\bbf-w||^{2n}}\cdot \sum_{s=1}^n(-1)^{s-1} (\overline{g}_s-\overline{w}_s) \overline{\partial g_1}\wedge\ldots\wedge\widehat{\overline{\partial g_s}}\wedge\ldots\wedge\overline{\partial g_n}.
\]
On the other hand, assuming ${\Gamma}_\varepsilon:=\{x\in B\,|\,|g(x)|=\varepsilon\}$ is smooth (which can always be achieved by varying $\varepsilon$), one has 
\begin{eqnarray}\label{w2}
(\ref{slr})=\frac1{(2\pi i)^n}\int_{{\Gamma}_\varepsilon} \frac{\omega}{(g_1-w_1)\ldots(g_n-w_n)}
\end{eqnarray}
where 
$
{\Gamma}_\varepsilon
$
is oriented by $d({\rm arg}\,\, g_1)\wedge\ldots\wedge d({\rm arg}\,\, g_n)>0$. (A proof of (\ref{w2}) for the regular values $w$ of $g$ can be found in \cite[Sect.5.17]{AGV}; by the principle of continuity for residues \cite[Ch.5, Sect.1]{GH} the formula holds true for all $w\in D_\varepsilon$.) To obtain (\ref{GGH}) it remains to apply the differential operator $\frac{\partial^{r_1+\ldots+r_n}}{\partial w_1^{r_1}\ldots\, \partial w_n^{r_n}}$ to the right-hand sides of (\ref{w1}) and (\ref{w2}) and set $w=0$.

\bigskip

\end{document}